\documentclass[a4paper,final,10pt]{article}  
\usepackage{times,amssymb,amsmath,graphicx,theorem,chicago}
\newcounter{propos}
\newcounter{proofs}
\theoremheaderfont{\scshape}
\theorembodyfont{\rmfamily}

\theoremstyle{break}
\newtheorem{propo}[propos]{Proposition}
\newtheorem{proof}[proofs]{Proof}

\begin{document}
\title{Viability Analysis of Management Frameworks for Fisheries} 
\author{K. Eisenack$^a$, J. Scheffran$^b$, and J.P. Kropp$^a$\\
\\
$^a$Potsdam Institute for Climate Impact Research \\
P.O. Box 60 12 03, 14412 Potsdam, Germany\\
eisenack@pik-potsdam.de, kropp@pik-potsdam.de\\
\\
$^b$University of Illinois, International Programs and Studies/ACDIS\\
505 Armory Ave, Champaign, IL 61820, USA\\
scheffra@uiuc.edu}
\date{}
\maketitle
\begin{abstract}
The pressure on marine renewable resources has rapidly 
increased over past decades. The resulting scarcity 
has led to a variety of different control and 
surveillance instruments. Often they have not improved the 
current situation, mainly due to institutional failure
and intrinsic uncertainties about the state of stocks.
This contribution presents an assessment of different management schemes
with respect to predefined constraints by utilizing viability theory.
Our analysis is based on a bio-economic model
which is examined as a dynamic control system in continuous time.
Feasible development paths are discussed in detail.
It is shown that participatory management may lead to serious problems if a purely resource-based management 
strategy is employed. The analysis suggests that a less risky management strategy can be implemented
if limited data are available.
\end{abstract}
\section{Introduction}
Marine fish stocks are under extreme pressure worldwide \shortcite{FAO.2001,Myers.2003}. Two divergent but closely related developments are observable \shortcite{Munro.1999,Pauly.2002}. On the one hand, increasing surveillance efforts, limited entries, or marine protected areas are established for mitigating overfishing, whilst on the other hand the fishing industry can be sustained at an economic level only
by high amounts of subsidies from the public sector \shortcite{Mace.1996,Banks.1999,Greboval.1999}. This is, in particular, remarkable since there has been an awareness of these problem for decades and most fisheries are
subject to management measures. The general situation in marine capture fisheries can be characterized by the following statements:
\begin{itemize}
\item Several authors argue that research as well as management instruments are overly biased 
toward the ecological viewpoint, while economic driving forces or the role of political decisions are 
rarely considered \shortcite{Lane.2000,Davis.2001}.
\item Fisheries are per se
associated with inherent uncertainty which is related to the partial opaqueness on both
the ecological and the economic systems \shortcite{Whitmarsh.2000}.
\item It is often mentioned that model approaches used to provide policy advice are  
based on unrealistic basic assumptions and/or inappropriate methodological concepts (cf. \shortciteNP{Imeson.2003}).
\end{itemize} 
However, one-sided views, distinct sources of vagueness, and unsuitable assessment techniques
not only restrict what can be modeled, but, they also impose severe limitations
regarding the explanatory power of the results. Thus, the basic question concerning the possibilities for a 
sustainable development in fisheries still remains: How can we manage fisheries in a way that avoids  ``hazardous developments''? Recently, so-called co-management schemes
have been introduced in order to mitigate the consequences of fishery mismanagement, but 
currently systematic assessments of control schemes regarding their effectiveness are
rarely employed. Considering this situation, it seems impossible to define safe management options for fisheries. 

One focal point of our examination is the inability to anticipate
exactly what will happen in the future of fisheries and marine stocks. Thus, we do not strive to determine
``optimal paths'' for the co-evolution of fisheries and marine resources, but rather for desirable
corridors\footnote{The conceptual ideas how these corridors can be used for policy advice and how judicious strategies for the management of development processes shall be designed, have been discussed in e.g.\shortciteNP{Kropp.1998a}. Inverse techniques (target-oriented approaches), such as
viability analysis, solve many analytical problems posed by the intricacy of the policy-bioeconomic complex.}. These corridors are constrained by measures representing
our knowledge of what should at least be avoided in order to achieve sustainability or prevent
catastrophic developments (for an example from climate research, cf. \shortciteNP{Petschel.1999}).  
We make use of these ideas by applying viability theory, which was developed 
by \shortciteN{Aubin.1991}, because it enables us to formally define the corridor boundaries. Such an analytical strategy allows to provide knowledge for decision-making,
although constraints are normative and an outcome of public discussions. 
In addition, it shifts
attention from the whole set of options to a constrained set of options, leaving  
space for adjustments acceptable in political frameworks.  
Additionally, we have combined the viability approach with game-theoretic
assumptions about the socio-economic mechanisms
constituting typical management frameworks. 
The management schemes we are analysing are an issue of the current debate
and are discussed briefly in the next paragraphs.

\begin{figure}[t]\centering
\includegraphics[width=1\linewidth]{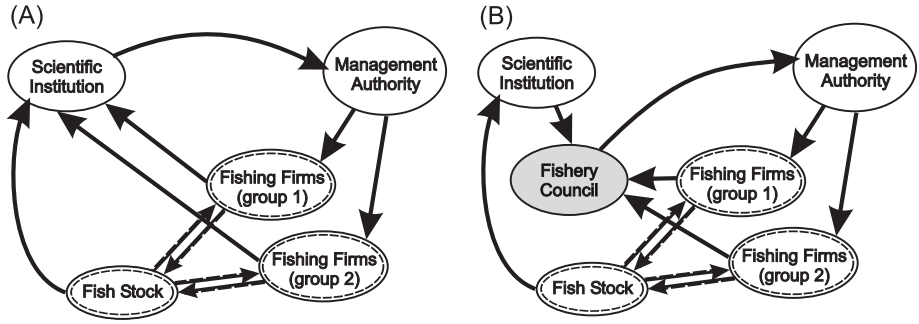}
\caption{Sketch of management schemes in fisheries:  (A) Top-down strategy: a scientific
organization collects data from fishing firms and own scientific trawls
and estimates the stock size. Resulting reports are delivered to the
management authority for decision support. The authority imposes
restrictions on the fishing firms. (B) For a co-management
a fishery council is introduced, in which representatives from science
and the fishing industry negotiate for appropriate catch
or effort restrictions. If they agree the resulting plan is
provided to the management authority, who approves it. 
An unmanaged fishery is represented by the
dashed cores. Here the firms calculate their efforts independently
and only related to the current stock status (measured by the catch from the previous
years, etc.).}
\label{fig.topdown}
\end{figure}

Without any regulation it is likely that in the future fish stocks
will be further depleted, as long as 
overexploitation is profitable for the individual fisherman. The revenues
from catches are private, while the costs induced by a reduced resource stock are
shared between all participants in the fishery (tragedy of the commons, cf. \shortciteNP{Hardin.1968}). 
Competing fishing firms base their decisions on deployable capital, necessary efforts, and 
on the observed state of the target species (see dashed core in Fig.\,\ref{fig.topdown}A,B). 
The latter is rather critical,
because the catch is commonly used as an estimator for the abundance of biomass.
This situation changes when a management authority introduces measures imposing
restrictions on the fishery (e.g. on gear type, allowable catches,
amount of effort, etc.). Such limitations change the decisions made by fishing firms 
(top-down management, cf. Fig.\,\ref{fig.topdown}A). Nevertheless, in many cases the results are disappointing, because 
restrictions are perceived as constraining economic opportunities.
Thus, fishing firms act as opponents to management authorities, often resulting in illegal landings
or mis-reported catches 
\shortcite{Agnew.2000,Hollup.2000,Robbins.2000,Jensen.2002}.  
Moreover, scientific
stock estimates are often not as reliable as required. Hence,
when a fishery reaches a state of crisis, scientific institutions come under
pressure in the public debate, i.e. for putting too much stress on
conservation objectives and neglecting economic sustainability. 

One solution to avoid these shortcomings is offered by co-management schemes, i.e. by including 
fishing firms in the decision-making process \shortcite{Jentoft.1998,Noble.2000,Charles.2001,Potter.2002}. 
If fishermen are involved in
the decision-making process, it is assumed that economic objectives
will complement conservational goals of governmental organizations.
It is the aim of this strategy that all actors in marine capture fisheries 
participate actively with respect to the overall target to keep the utilization of
marine resources sustainable. Self-governance within
a legal framework constituted by governments is a basic principle of this strategy.
As a result, the fishing industry
represents economic objectives in a negotiation process with agents pursuing
conservational goals.  Before such a situation comes into
play, a variety of conflicts often have to be resolved, e.g. conflicts over who
owns and controls access, how policy and control mechanisms should be 
carried out, or conflicts between local stakeholders  (for a 
nomenclature of common conflicts in fisheries, cf. \shortciteNP{Charles.1992}).
Typically, this type of management is exercised via a fishery council 
where the representatives of fishing firms, processing firms, scientific
institutions and policy negotiate, e.g. about the total allowable
catch for each different species. This plan has to be approved by a
governmental authority and is executed by a management organization which operates
in close collaboration with local fishermen (see Figure \ref{fig.topdown}B).
A plausible conjecture is that under co-management fishermen will
show higher compliance with the resulting constraints (cf. \shortciteNP{Pinkerton.1989};
\shortciteNP{Mahon.2003}). 

However, an additional problem occurs if management regimes consider only 
biological measures as steering targets. Such a strategy has been claimed as
'ichthyocentrism', indicating that scientific advice puts too much emphasis on the
resource itself (the various fish stocks, in particular their biomass) \shortcite{Lane.2000,Davis.2001} 
compared with efforts to examine the behaviour of the resource users, their economic settings,
and aims.
In such a case, the management success cannot depend on an exact stock assessment, 
which is impossible, due to insufficient measurement and sampling methods.  

Our study applies viability analysis in order
to assess different management schemes in marine fisheries. Recent studies
have shown that this methodology is valuable for determining completely unacceptable outcomes and defining judicious measures for mitigation
\shortcite{Petschel.1999,Moldenhauer.1999,Bene.2001,Aubin.2005}.
Therefore, we try to find out how management in fisheries 
should be designed under defined constraints in order
to achieve safe limits.

Pursuing this goal we have defined a model, including stock dynamics as well as 
economic and political decision-making.
We motivate the viability constraints and use the catch recommendations of the scientific
institution participating in the co-management process as control variable.
The viability criteria are imposed on the model
in order to assess and develop different control strategies.  
A discussion of the results and a summary concludes the paper.

\section{The Dynamic Model}
The basic state variable of the model is the biomass of a
fish stock $x$, which is influenced by the total
harvest $h$ in the fishery. If we introduce a
recruitment function $R$ which assigns the growth of the biomass
to a given stock, we obtain
\begin{equation*}
  \label{eq:ode_stockdyn}
  \dot{x} = R(x) - h
\end{equation*}
as ODE for the stock dynamics.
As usual, $R$ is assumed to be of the Schaefer type
\shortcite{Schaefer.1954}, producing logistic growth
with maximum $\bar{R}$ at $x_{MSY}$ and $R(0)=0$.
\begin{equation*}
\forall x \in [0, x_{MSY}]: R(x) > 0 \quad\text{and} \quad D_x R(x) > 0,
\end{equation*}
where $D_x$ denotes the differential operator with respect to
the argument $x$. For $x > x_{MSY}$ we have $D_x R(x) <0$.
Due to the complexity of ecosystems (embedding of target species
into several food webs, large scale impacts such
as ENSO, etc.) we have to deal with limited knowledge about the behaviour of a fish stock. 
Thus, no additional assumptions about $R$ are introduced.

The next step is to determine the amount of total harvest $h$. In the presented model
we concentrate on output-management, i.e.  the negotiation process is
about the allocation of catch quotas $q_i$ to groups of fishing firms
$i=1,\ldots,n$.  We assume that it is profitable for each firm to realize a
harvest which it is allowed to catch.  Thus, the resulting total harvest is $h=\sum_i
q_i$.
To model the negotiation process a game theoretic approach can be applied where
the fishing actors agree on particular limits for harvest and the allocation of
efforts (cf. \shortciteNP{Scheffran.2000}). Expanding this
approach, a scientific institution and representatives from the fishing
industry ``bargain'' for the total harvest $h$ and the individual quotas $q_i$.
When these pressure groups agree on an allocation, the result is transformed into practice by
the management authority.  It is further presumed that the negotiations are opened by the
scientific institution, which makes a recommendation $r \ge 0$ for the total
catch.  Each group of the fishing industry tries both (i) to
get a share as high as possible of the total harvest $h$ and (ii) to increase $h$ above the catch
recommendation $r$ in order to improve their profits. 

The optimal quota and optimal increase $h-r$ may differ among the
groups, e.g.  due to their capitalization or technical efficiency
(artisanal vs. industrial fishery). But, there is also a 
trade-off between higher profits resulting from higher quotas and
deviation costs $d_i$ imposed by exceeding the scientific
recommendation. These costs are linked to the legitimation of
bargaining positions challenging the scientific advice and increasing
transaction costs of fierce negotiations (time, expertise, human and
social capital, data retrieval, public relations, etc.). How strong
this trade-off is, depends i.a. on reputation, the political influence
and the availability of information -- all of which may differ between
the pressure groups. We further assume that the fishing groups act
myopically, which is realistic if the influence of single fishing
firms on the resource is neglectable \shortcite{Banks.1999,Kropp.2004d} and
if they push their representatives for higher quotas. This entails
that they only account for short-term deviation costs. Myopic
behaviour can also be observed if the burden of long-term
responsibility is shifted to the scientific institution (allowing
fishing firms to reduce subjective uncertainties). 

We denote the profit of a group $i$ as $\pi_i$. It depends on
the quota $q_i$ and the available amount of fish $x$ (which is the
same for all fishing firms). This function also represents the
efficiency of boats, fishing gear, technological equipment,
etc. In a situation without deviation costs,
the function has the form
\begin{equation*}
  \pi_i(q_i,x) =  p \, q_i - c_i(q_i,x).
\end{equation*}
The first term represents revenues on markets, where $p$ corresponds
to the market price (which is assumed to be exogenous), while $c_i$ is
a cost function, assigning variable costs to a realized harvest $q_i$. It is
economically reasonable to assume that $c_i$ is increasing in $q_i$,
since more labour, time, fuel etc. is needed for a larger catch. On
the other hand, costs are decreasing in $x$ due to higher densities of
fish.
Each firm individually selects 
catch $q_i$ to obtain an optimal profit for a given price and target
species. However, in a co-management framework the profit function is
modified by deviation costs:
\begin{equation*}
  \pi_i(q_i,x) = p q_i - c_i(q_i,x) - d_i \bigg(\sum_{j=1}^n q_j - r \bigg).
\end{equation*}
It should be noted that the deviation costs not only depend on the
individual decision $q_i$, but also on the quota allocated to the other
pressure groups. If $\sum_j q_j - r$ becomes negative, we assume
that $d_i$ vanishes, since deviation costs do not apply if the sum of all
quotas is below the scientific recommendation. It is reasonable to
presume that each
$d_i$ is a monotonic increasing function.

Now, the Nash equilibrium of the negotiation process is given by a
quota allocation which assigns an individual quota $q_i$ to each group $i$ 
that maximizes $\pi_i$ with respect to $q_i$ for given $p$, $x$ and
the quotas of the other participants. The resulting total harvest is the sum of
all quotas.
By applying basic calculus it is shown 
that all $\pi_i$ are concave and
continuously differentiable with respect to $q_i$. Solutions of the equation
system 
\begin{equation*}
\forall i=1,\ldots,n: D_{q_i} \pi_i = 0
\end{equation*}
are such negotiation equilibria.

In the following analysis we restrict this general approach to the
case of two specific fishery groups and provide a possible functional
specification for variable and deviation costs:
\begin{equation*}
  c_i(q_i, x) := \frac{\alpha_i q_i + \beta_i q_i^2}{x}\ .
\end{equation*}
\begin{equation*}
  d_i(q_1 + q_2 - r) :=
        \begin{cases}
        0 \quad \text{if} \quad q_1+q_2 < r \\
        \kappa_i(q_1+q_2-r)^2 \quad \text{otherwise} \ .
        \end{cases}
\end{equation*}
The parameters $\alpha_i, \beta_i, \kappa_i$ $(i=1, 2)$ are not completely known, but positive.

First assume that we are in the case of non-vanishing deviation costs.
If $x$ is positive and since $q_i \ge 0$ the resulting profit
functions $\pi_i$ 
are continuously differentiable and concave with respect to $q_1, q_2$. Thus, the
Nash equilibrium (for given $p, x, r$) is obtained by solving
\begin{align*}
D_{q_1} \pi_1 = p - \frac{\alpha_1 + 2 \beta_1 q_1}{x} - 2 \kappa_1
(q_1+q_2-r)         &= 0, \\
D_{q_2} \pi_2 = p - \frac{\alpha_2 + 2 \beta_2 q_2}{x} - 2 \kappa_2
(q_1+q_2-r)         &= 0, \\
\end{align*}
for $q_1, q_2$. The latter equations define a total (binding)\footnote{The two cases
considered here comprise binding and non-binding recommendations. The former indicates
that the achieved harvest is exactly constrained by the recommendation, the latter that the
harvest is lower than the recommendation, i.e. the industry catches -- voluntarily --
less than suggested.} harvest
\begin{equation}
 h_b(x,r) = q_1 + q_2
 = \frac{ u p x + w x r - v}{\beta_1\beta_2 + w x},
\label{eq:h1}
\end{equation}
where
\begin{align*}
  u &:= \frac{1}{2} (\beta_2+\beta_1) > 0, \\
  w &:= \beta_1\kappa_2 + \beta_2\kappa_1 > 0, \\
  v &:= \frac{1}{2}(\alpha_1\beta_2+\alpha_2\beta_1) > 0.
\end{align*}
\begin{figure}[ht]\centering
\includegraphics[width=.8\linewidth]{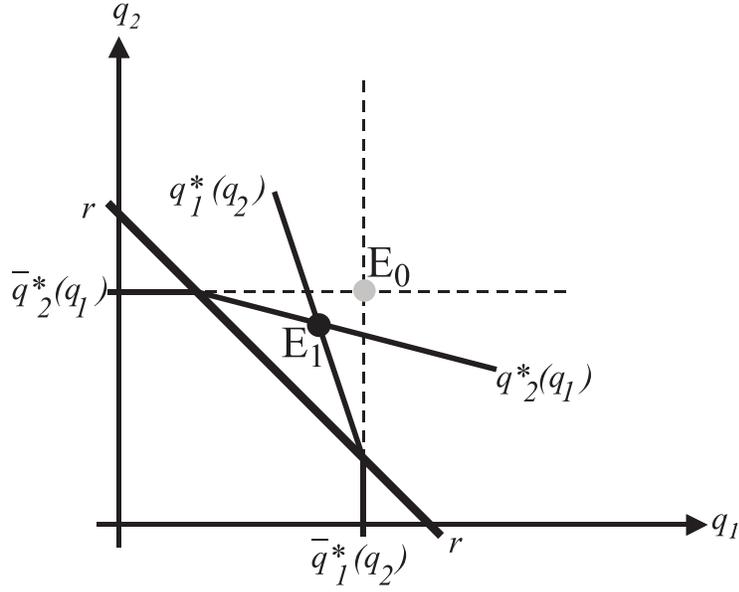}
\caption{Negotiation equilibrium of the quota setting game. $q^*_1$
and $q^*_2$ denote the optimal catch for both groups if there were no
management (i.e. deviation costs $d_i \equiv 0$). In this case
the choice of the groups is independent from each other (indicated
by dashed lines), and will
result in the equilibrium $E_0$. By introducing a catch
recommendation $r$ all points to the upper right of the diagonal are
associated with positive deviation costs. Therefore, the optimal
choice of each group depends on the choice of the other (indicated by
solid lines) and results in the equilibrium $E_1$.}
\label{fig.equilib}
\end{figure}
The binding harvest varies according to
\begin{equation}\label{eq:Drh1}
D_r h_b =  \frac{w x}{\beta_1\beta_2 + w x} > 0,
\end{equation}
and
\begin{equation}\label{eq:Dxh1}
D_x h_b =
 \frac{vw + (up+wr)\beta_1\beta_2}{(wx+\beta_1\beta_2)^2}
 > 0 \ .
\end{equation}
If the catch recommendation $r$ is so high 
that it is not profitable for fishing firms to exceed it, then no
deviation costs apply in the Nash equilibrium. This special case $\kappa_i = 0$ 
results in a total harvest (associated with a non-binding recommendation)
\begin{equation}
 h_n(x,r) 
 = \frac{ u p x - v}{\beta_1\beta_2} =: \hat{r}(x).
\label{eq:h2}
\end{equation}
Thus, for $r \ge \hat{r}(x)$, we have $h=h_n(x,r)$, i.e. the industry catches voluntarily less fish than recommended and therefore, no deviation costs apply. For $r < \hat{r}(x)$ the harvest
is $h=h_b(x,r)$.
Summarizing we obtain the total harvest function\footnote{There are two additional specific 
cases that have to be considered (cf. Eq.\,\protect\ref{eq:harvest}: $h(x,r)=0$). 
If $\hat{r}(x) \le 0$, it is not profitable to cast for fish, even if there were no
deviation costs, and therefore $h=0$. If
$0 \le r \le \hat{r}(x)$, but $h_b(x,r) \le 0$,
the recommendation is so tight that it prevents commercial fishing activities.}
\begin{equation}
  \label{eq:harvest}
  h(x,r) =
\begin{cases}
h_b(x,r) & \text{if}\quad r \le \hat{r}(x) \quad\text{and}\quad h_b(x,r) \ge
0, \\
h_n(x,r)=\hat{r}(x) & \text{if}\quad r \ge \hat{r}(x) \ge 0, \\
0  & \text{if}\quad \hat{r}(x) \le 0, \\
0  & \text{if}\quad 0 \le r \le \hat{r}(x)  \quad\text{and}\quad  h_b(x,r) \le 0.
\end{cases}
\end{equation}
We remark that
\begin{equation}
  \label{eq:hath1h2}
  h_n(x,r)=  h_n(x,\hat{r}(x)) = \hat{r}(x) = h_b(x,\hat{r}(x)).
\end{equation}
As $h_b$ is monotonically increasing in $r$ while $h_n$ is
independent from recommendations, we find
\begin{align}\label{eq:h2geh1}
  r > \hat{r}(x) & \Rightarrow h_n(x,r) < h_b(x,r), \\
  r < \hat{r}(x) & \Rightarrow h_b(x,r) < h_n(x,r),
\end{align}
and conclude that always
\begin{equation}
  \label{eq:hlehatr}
  h(x,r) \le \hat{r}(x).
\end{equation}
Equality is only possible if non-binding
recommendations are made ($r \ge \hat{r}(x)$).
Otherwise, the result of the negotiation process
is below the harvest which would be economically
optimal in the case of absent deviation costs (cf. Fig. \ref{fig.equilib}).
Since additionally 
\begin{equation}
  \label{eq:Dxrhat}
  D_x h_n = D_x \hat{r} = \frac{up}{\beta_1 \beta_2} > 0,
\end{equation}
the total harvest $h$
is increasing (not strictly) with an increased abundance of fish (supposing the
recommendations are unchanged).
\section{Viability Constraints for Sustainability}
To answer the question whether a fishery described by this model
can be managed in a sustainable way or not, it is necessary to
specify this objective in more detail. Generally, sustainability can be
characterized by ecological, economic, and
social dimensions. Here, we concentrate on the first two and facilitate their formalization 
in the framework of viability theory \shortcite{Aubin.1991}. 
Viability constraints characterize an acceptable sub-region of the phase space.
A time evolution of a system is called viable (or
sustainable) if it remains in this region indefinitely.
If a development process is controlled, in the examined case by the harvest
recommendation $r$, we want to analyse whether a control strategy keeps it
viable or not.

The choice of viability constraints cannot be purely justified by
empirical considerations, because it involves value-laden normative
settings (e.g. on what is at least acceptable or what do we want). The viability
concept allows evaluation of different normative settings with 
respect to their consistency and consequences, for instance, whether a given management
framework admit controls which satisfy the constraints.
For our examination of marine fisheries two reasonable viability constraints 
are defined and investigated. We deduce 
conditions under which a control rule for $r$ exists, respecting
both constraints at the same time:
\begin{enumerate}
\item Ensure that the biomass of a stock resides always above a minimal
level $\underline{x} > 0$, i.e.
\begin{equation*}
  \forall t: x(t) \ge \underline{x}.
\end{equation*}
\item  Require that a minimum total harvest $\underline{h} > 0$ can always
be realized or exceeded, i.e 
\begin{equation*}
  \forall t: h(t) \ge \underline{h}.
\end{equation*}
This
harvest covers fixed costs in the fishery,
guarantees a minimum level of employment, or sustains food safety.
\end{enumerate}
In the following, we refer to the first criterion as ``ecological'' and
to the second one as ``economic'' viability. The viability problem is
introduced more formally as follows.

Let $F(x)$ denote the set of all derivatives $\dot{x}$ which are
admissible in $x$, e.g. $F(x) := \{ R(x) - h(x,r) \mid r \ge 0 \,\text{and}\, h(x,r) \ge \underline{h} \}$ 
for all changes of fish stocks resulting
from an economically viable harvest recommendation $r$. An interval
$I=[\underline{x}, \infty]$ is called a viability domain of $F$ if 
\begin{equation}
  \label{eq:viabdomain}
\forall x \ge \underline{x}: F(x) \neq
\varnothing \,\text{and}\,
 F(\underline{x}) \cap [0, \infty] \neq \varnothing.
\end{equation}
The viability theorem \shortcite[p.91]{Aubin.1991} states that if $I$ is a
viability domain of $F$, 
then for every initial value $x_0 \in I$ there exists a control path
$r(\cdot)$, 
such that the solution of the initial value problem
\begin{align*}
  \dot{x}(t) &= R(x(t))-h(x(t),r(t)) \in F(x), \\
   x(0) &=x_0,
\end{align*}
remains in $I$, i.e. $\forall t\ge 0: x(t) \in I$
(more technical premises of this theorem
can be found in \shortciteN{Aubin.1991}).

In the following we determine criteria for $I$ being a viability domain,
depending on given values of $\underline{x}$ and $\underline{h}$.
For such values
a ``wise'' harvest recommendation can keep the fishery within economic
and ecological sustainable limits.
If the harvest recommendation $r$ is binding then $\dot{x} \ge 0$, if $R(x)
\ge h_b(x,r)$ (catch matches recruitment). The
solution of this inequality yields
\begin{equation}
 \label{eq:vaibdotx}
  r \le \frac{ R(x) (wx+\beta_1\beta_2) + v} {wx}
          - \frac{up}{w} =: \bar{r}(x),
\end{equation}
i.e. a maximal recommendation which results in a negotiation
equilibrium for the total harvest which is low enough to ensure an increasing
fish stock. This, however, depends on the state of the resource, the market, and
other economic parameters.
Relation (\ref{eq:vaibdotx}) does not necessarily have
a positive right-hand side. Since $r \ge 0$, it is impossible to
stabilize the fish stock in such a case. 
The
maximal recommendation does not necessarily increase with the fish stock.
Differentiating with respect to $x$ yields
\begin{equation}\label{eq:Dxrbar}
D_x \bar{r} = \frac{w x^2 D_x R(x) + \beta_1\beta_2 (x D_x R(x) - R(x)) - v}{w x^2},
\end{equation}
where the expression
$x D_x R(x) - R(x)$ may become negative.

In the case of non-binding recommendations (catch is lower than recommendation) $r \ge \hat{r}(x)$, we
need $h_n(x,r)= \hat{r}(x) \le R(x)$ for an increasing resource. This condition
does not depend on the concrete value of $r$. 

Due to relation (\ref{eq:viabdomain}) all of the above conditions only have to hold at $x=\underline{x}$,
and we can determine the maximal harvest recommendation if
$R(\underline{x})$ is known. Again, it may happen that, e.g. due to high
prices, this condition cannot be fulfilled. It depends on a non-trivial
relation of economic and ecological parameters, whether the quota
negotiation process can yield a viable result.
To satisfy the economic viability constraint, $r$ has to be chosen for every time such that $h(x,r)\ge \underline{h}$.
For $r \le \hat{r}(x)$ we solve $h_b(x,r)\ge \underline{h}$ for $r$, yielding
\begin{equation}
  \label{eq:viabh}
  r \ge \frac{ \underline{h} (w x +\beta_1\beta_2) + v}
               {w x}
          - \frac{up}{w}  =: \underline{r}(x).
\end{equation}
Contrary to the former case,
we obtain a lower
limit for the harvest recommendation. While the similarity of the right hand
sides of Eq. (\ref{eq:vaibdotx}) and Eq. (\ref{eq:viabh}) is obvious,
the partial derivative with respect to $x$ simplifies to
\begin{equation}\label{eq:Dxru}
 D_x \underline{r} = - \frac{\underline{h}\beta_1\beta_2+v}{wx^2} < 0,
\end{equation}
i.e. for larger fish stocks, lower catch recommendations guarantee
economic viability. However, for non-binding recommendations, the
condition $\hat{r}(x) \ge \underline{h}$ has to hold for every $x\in
I$.

Combining both viability constraints, we must distinguish whether
there is the possibility to choose $r$ at $x=\underline{x}$ 
such that it respects
Eq. (\ref{eq:vaibdotx}) and Eq. (\ref{eq:viabh}) (in the case of
binding recommendations), i.e.
$\underline{r}(\underline{x}) \le r \le \bar{r}(\underline{x})$,
which is equivalent to $\underline{h} \le R(\underline{x})$.
Summarizing all possible cases yields (cf. appendix \ref{prop:vd} for the proof):
\newline
\newline\noindent The~interval $I=[\underline{x},\infty]$ is a viability domain, if and only if
\begin{alignat}{3}
\label{prop:vda}
 & \text{(i)}&\qquad &\hat{r}(\underline{x}) \ge \underline{r}(\underline{x}), \quad\text{and} \\  \nonumber
 & \text{(ii)}&\qquad &\underline{h} \le R(\underline{x}), \quad\text{and}\\
 & \text{(iii)}&\qquad &\underline{r}(\underline{x}) \ge 0\quad\text{or}\quad \bar{r}(\underline{x}) \ge 0.  \nonumber
\end{alignat}
Although for the proof (\ref{proof:vd}) an awkward set of cases has to be considered, these conditions are easy to interpret.
The compatibility of both viability criteria only depends on
the relation between the recruitment function, the aspiration level
for harvest and the efficiency of the firms.
If the recruitment at the minimal viable stock
level $\underline{x}$ lies below the required harvest, there exists an obvious
contradiction
between economic and ecological targets (condition \ref{prop:vda}ii).
In this situation,
stock approaches $\underline{x}$, the fisheries council has to decide
whether to sacrifice conservational or harvest objectives.
According to (\ref{prop:vda}i) it is profitable to fish at least
 $\underline{h}$, even if no deviation costs apply. 
Otherwise no harvest suggestion guarantees an adequate yield,
because the harvest is always below or equal to $\hat{r}(x)$.
Condition (\ref{prop:vda}iii) can be met, if the recommendations can reduce the
harvest below $R(\underline{x})$ in case the stock approaches
$\underline{x}$. Otherwise, a recommendation $r < 0$ would be needed
in order to achieve a sufficient harvest reduction $h(x,r)< R(x)$. 
However, in our examinations this makes no sense.
\section{Assessment of Recommendation Strategies}
As long as the system stays in the viability domain it is {\em
possible} to keep it viable if an appropriate strategy is
selected. However, this does not ensure that {\em every} control
strategy is successful. In addition, in a critical state in which the viability
constraints are not satisfied, it is not clear 
whether a selected strategy leads to ill-management or 
forces a fishery back to sustainable limits.
Consistently, we examine and assess the viability of different
recommendation strategies in this section. 
Formally, such a strategy assigns a value for $r$ to a given system
state, i.e. a closed-loop control according to the following schemes.
\begin{itemize}
\item {\em Ichthyocentric control}: The harvest recommendation is purely
based on an estimate of the stock recruitment.
\item {\em Conservative control}: The harvest recommendation is based
on economic viability in the sense that the minimum harvest
$\underline{h}$ is always realized.
\item {\em Qualitative control}: In this case, due to uncertainties,
recommendations are only
based on qualitative economic and ecological observations.
\end{itemize}
\subsection{Ichthyocentric Control}
The first control strategy is called ``ichthyocentric'', in order to make clear
that the scientific organization considers only the state
of a fish stock for their harvest suggestions, but not the socio-economic conditions.
It is assumed that the harvest recommendation $r$ equals
the estimated recruitment of the targeted species and 
that the scientific institution is able to  
estimate correctly the recruitment function, i.e.
\begin{equation*}
 r=R(x).
\end{equation*}
This is a challenging task, since an exact estimation of the stock
biomass is bound to fail, due to unavoidable measurement deficits.
However, let us assume that
the estimator is roughly correct.
The examination of the ichthyocentric control strategy 
shows that even in this ideal case economic
viability is not always guaranteed. 
If $r= R(x) < \underline{r}(x)$ the recommendations are too low to
sustain a harvest rate $\underline{h}$
(in particular, if the fishery stays in a viability domain, cf. Prop. \ref{propo:icht1} in the appendix).
The situation is even more dramatic if one focuses on ecological viability.
If recommendations $r$ are binding our results show that the
catch exceeds $r$. As a consequence, recommending $R(x)$ leads to a
catch above recruitment which implies an always decreasing fish stock 
($h(x,R(x)) \le R(x)$, if and only if $r=R(x) \ge \hat{r}(x)$, cf. Prop. \ref{propo:icht2}).

Ichthyocentric control is only sustainable
if recruitment is always large enough to allow a minimal harvest and
if recommendations are not binding at $x=\underline{x}$.
The latter means that the realized catch $\hat{r}$ must be significantly lower than
the scientific recommendation, a situation which 
normally does not occur in industrial capture fisheries.
But it might be observable in low capitalized fisheries, e.g. if the competitors
only use small boats.
Otherwise, the stock will necessarily decrease below
$\underline{x}$. If this is already the case, the chance for a
regeneration of stocks is rather bad until $\hat{r}(x)<R(x)$ is approached.
On this level the firms voluntarily catch less than recruitment.
We can summarize that even in the case of a perfect stock assessment, the
ichthyocentric strategy exposes the fishery to a risky development
path. 
\subsection{Conservative Control}
\begin{figure}[t]\centering
\includegraphics[width=.495\linewidth]{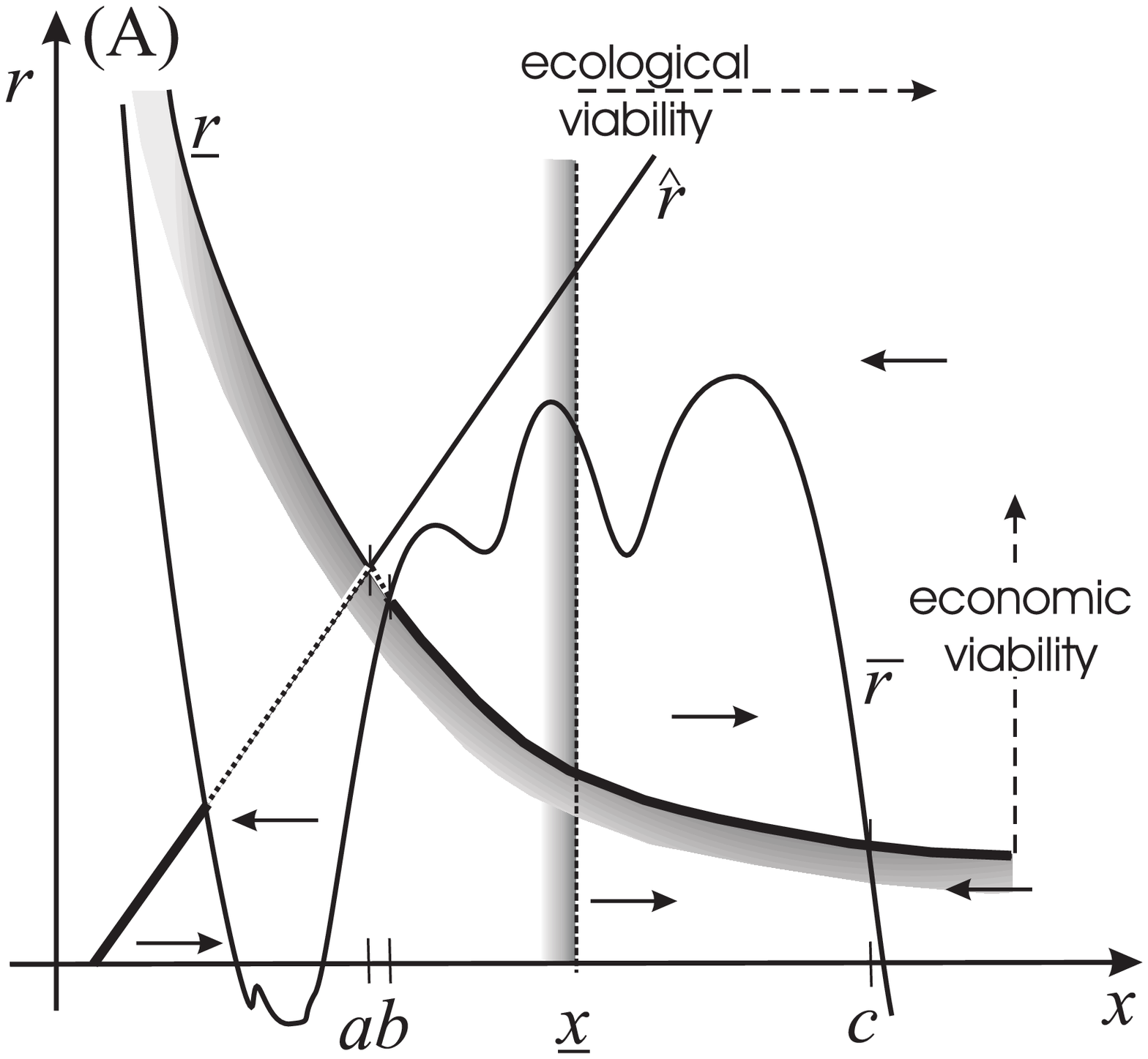}
\includegraphics[width=.495\linewidth]{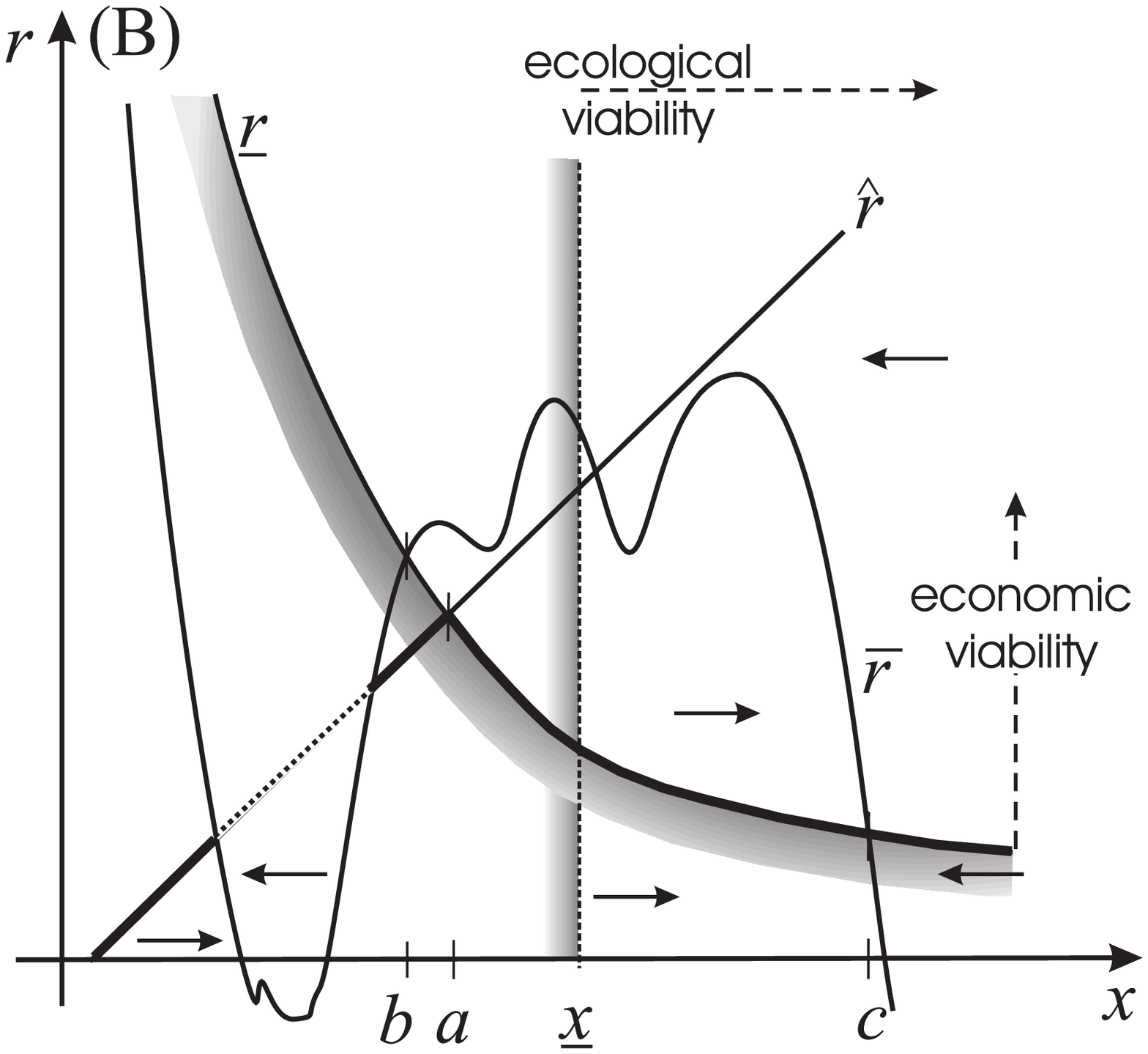}
\caption{
Phase portraits for the conservative control strategy: (A) $a<b$ left, (B) $b<a$. 
Below $\underline{r}$ a viable harvest cannot be obtained. Solid arrows denote the
direction of $x$, being positive below $\bar{r}$ and negative
above. Note that it is possible for both cases that $\bar{r}$ has no zero (then
the separated areas below $\underline{r}$ are joined) or multiple zeros. In addition,
if $a<b$ it might be possible that $\bar{r}$ becomes larger than $\hat{r}$ (cf. (A)), then
the area above $\bar{r}$ decomposes into multiple sections which is not shown here, because
it does not change the general statements). The bold and
dotted lines correspond to recommending $\underline{r}$, resulting in
an increasing (straight) or decreasing stock (dotted). Note that
for $x<a$, recommending $r=\underline{r}(x)$ results in the same catch as recommending $r=\hat{r}(x)$.
}
\label{fig:phase12}
\end{figure}

The conservative control strategy aims to ensure that the economic
viability criterion is satisfied, but nothing more: the scientific
institution always approves
\begin{equation*}
  r=\underline{r}(x).
\end{equation*}
To evaluate this strategy the phase space
structure of the model must be discussed in more detail
(Fig. \ref{fig:phase12}). 
We define the fish stock level $a$ by the intersection 
$\underline{r}(x)=\hat{r}(x)$. It is unique since
$\hat{r}(x)$ is linearly increasing, while 
$\underline{r}(x)$ is a monotonically decreasing function
(Eq. \ref{eq:Dxru}).
If $[\underline{x}, \infty]$ is a viability domain and
$\hat{r}(\underline{x}) \ge \underline{r}(\underline{x}) \ge 0$
(Eq. \ref{prop:vda}) holds, it follows that $a < \underline{x}$ and
$\hat{r}(a)=\underline{r}(a)>0$. 
We further define the fish stock levels $b < c$ as the solutions of
$\bar{r}(x)=\underline{r}(x)$. Although $\bar{r}(x)$ is (in general)
not monotonic, this equality only holds if
$\underline{h}=R(x)$  (cf. Eq. \ref{eq:vaibdotx} and \ref{eq:viabh}). 
Since the recruitment function $R$ is of Schaefer type, only a
maximum of two solutions can occur.
If $[\underline{x}, \infty]$ is a viability domain, with the consequence that $R(\underline{x}) \ge \underline{h}$ holds (Eq. \ref{prop:vda}),
these two intersections satisfy $b < \underline{x} < c$
(the only exception is the marginal case $R(\underline{x})=\underline{h}$).
We have, in principle, to distinguish the cases $a < b$ and $b < a$ 
(cf. Fig. \ref{fig:phase12}A,B). We further observe
that $\lim_{x \rightarrow 0} \bar{r}(x) = +\infty$, indicating
an additional intersection between $\hat{r}$ and $\bar{r}$ below $b$. 
For the conservative control scheme the system always evolves along the graph of the
bold straight and dotted line (cf. Fig. \ref{fig:phase12}).

If the fishery stays within the viability domain, we always have
$r=\underline{r}(x) \le \hat{r}(x)$
(Prop. \ref{prop:vd}) and
$h_b(x,\underline{r}(x))=\underline{h}>0$ holds. Consequently, 
it follows from Eq. (\ref{eq:harvest})
that $h(x,\underline{r}(x))=h_b(x,\underline{r}(x))\ge\underline{h}$,
i.e. conservative control guarantees economic viability in the
viability domain. Also, due to the fact that within a viability domain 
$h(\underline{x},\underline{r}(\underline{x}))=\underline{h} \le
R(\underline{x})$ is valid, the ecological viability constraint is met.
Therefore, it can be deduced that conservative control is -- in general -- a viable
strategy. Nevertheless, the
following two aspects have to be considered. (i) This strategy only
guarantees a minimal catch level, although a larger harvest might be
viable, too (for $x \in [\underline{x},c]$). (ii) Taking into account the deficits 
in exact stock estimation, it may happen that the state of the stock is
already below $\underline{x}$ indicating that no viable control
exists anymore, but this fact is unrecognized by the resource users.
However, it is of specific interest whether a stock will recover to a viable level or not
if the conservative control strategy is still exercised in such a situation. 
Both cases are possible, in some cases it the stock
can recover, but also worser scenarios can occur, e.g.
if $\underline{r}(x) \le \hat{r}(x)$ (as in the viability domain), but
$R(x) < \underline{h} \Leftrightarrow \bar{r}(x) < \underline{r}(x)$
(cf. Fig. \ref{fig:phase12}A, dotted line), i.e. recruitment is below 
minimal harvest and the stock further decreases.
Another situation is given if $\bar{r}(x) < \hat{r}(x) < \underline{r}(x)$
(cf. Fig. \ref{fig:phase12}B, dotted line), where
$h(x,\underline{r}(x))=\hat{r}(x)=h_b(x,\hat{r}(x)) >
h_b(x,\bar{r}(x)) = R(x)$.

We conclude that conservative control can satisfy economic and
ecological viability criteria, but only for non-critical situations (i.e.
in the viability domain). 
In contrast to the ichthyocentric control scheme, it has the advantage
that only qualitative information is needed to exercise this type of
management, i.e. the scientific institution
only has to supervise whether the realized catch of the fishery is above or
below $\underline{h}$. In the former case, catch recommendations
should be reduced, while in the latter they should be increased. 
This perspective considers the problems 
arising from uncertainty in fisheries management more seriously.
However, in a crisis (i.e. being outside the viability domain) 
where the aspiration level for harvest is too high or the abundance
of the targeted species is rather small for a viable control, we cannot be 
assured that the resource recovers by applying this management strategy.
But even in the viable case, profits in the fishery are still limited
to a minimum. 
\subsection{Qualitative Control}
The third alternative extends the qualitative view discussed above in
order to increase profits and limit risks in critical
situations, i.e outside the viability domain.
It is based only on qualitative observations, which means
that the exact numerical 
values of $x, R, \underline{r}$ and $\bar{r}$
are not known, but that 
it can be determined correctly 
whether $x$ is decreasing or increasing in
time, whether the realized catches exceed $\underline{h}$ or
not, and whether recommendations are binding or not. 
The applied control rule provides only qualitative advice, i.e.
whether $r$ has to be increased or decreased. The proposed 
qualitative control strategies are summarized in Tab. \ref{tab2}.
\begin{table}[h]\centering
\caption{Qualitative control as discussed in detail in the text. The
  numbers of rules corresponds to areas in the phase portrait
  (cf. Fig. \protect\ref{fig:phase12a}). Note that the coat indicates
the history of a fishery (cf. text). In the phase plots they are indistinguishable.}
\label{tab2}
\begin{tabular}{lll} 
\hline
rule no. & qualitative observation & reaction \\
\hline
(0) & $r>\hat{r}(x)$ &  decrease $r$ \\
(1) & $x$ increases and $h > \underline{h}$ & increase $r$ \\
(2) & $x$ increases and $h < \underline{h}$ & increase $r$ \\

(3') & $x$ decreases and $h > \underline{h}$ and mature & decrease $r$ \\
(4') & $x$ decreases and $h < \underline{h}$ and mature & moratorium: $h=0$ \\

(3'') & $x$ decreases and $h > \underline{h}$ and emerging & decrease $r$ until $h = \underline{h}$ \\
(4'') & $x$ decreases and $h < \underline{h}$ and emerging & increase $r$ until $h = \underline{h}$ \\

\hline
\end{tabular}
\end{table}
\begin{figure}[t]\centering
\includegraphics[width=.495\linewidth]{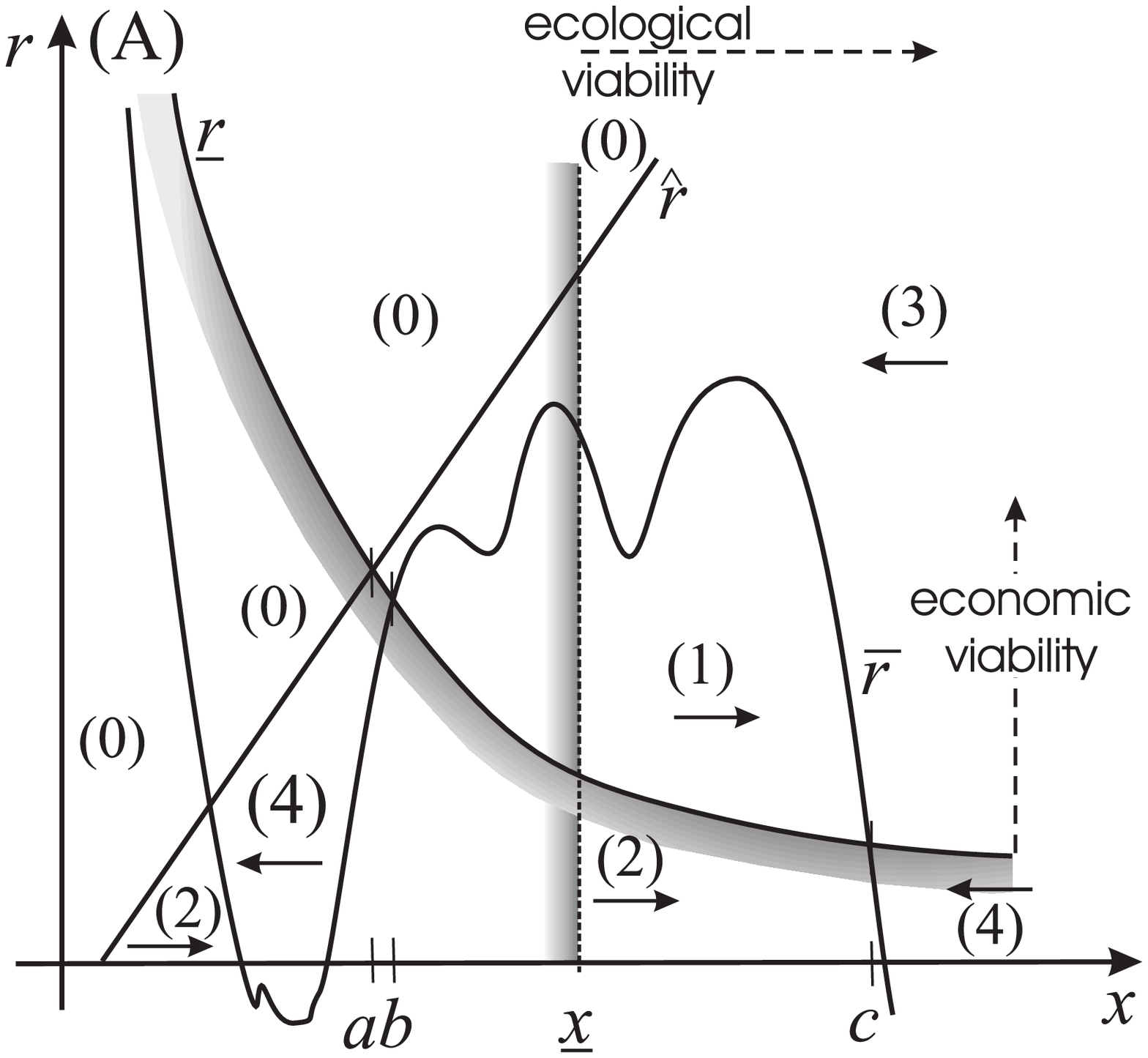}
\includegraphics[width=.495\linewidth]{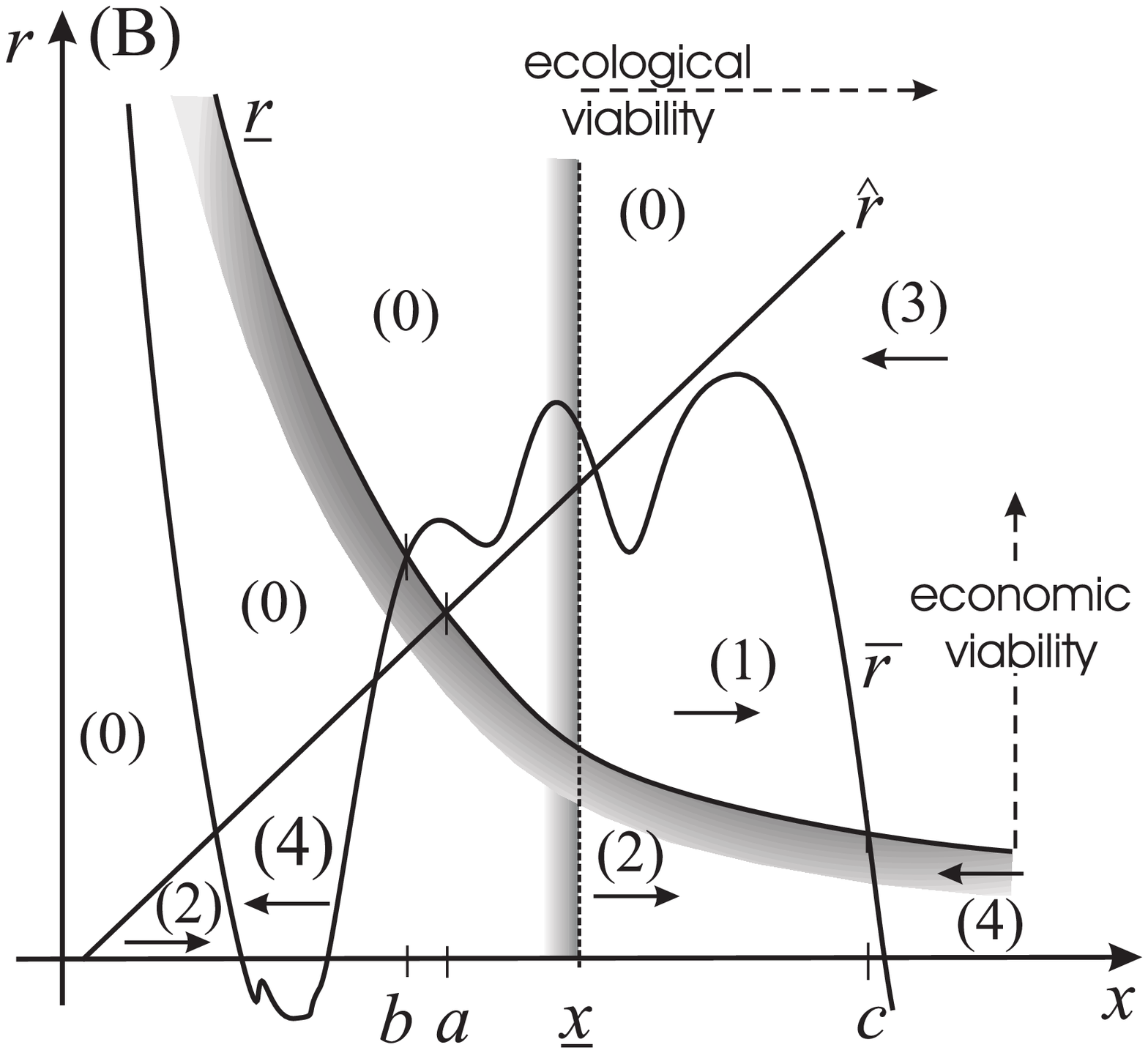}
\caption{
Phase portraits for the qualitative control strategy: (A) $a<b$, (B) $b<a$. 
The numbered areas correspond to situations which can be 
distinguished qualitatively and are discussed in detail in the text.
}
\label{fig:phase12a}
\end{figure}
To discuss this control scheme in detail, we refer to phase space portraits again
(Fig. \ref{fig:phase12a}). The underlying idea  is
to combine qualitative observations regarding a fishery with phase
plane analysis. This knowledge makes it feasible to decide
in which region of the phase plane a fishery is currently situated.
If sufficient knowledge is available, we additionally use information
on the history of a fishery in order to decide this question. This makes it feasible
to classify it either as emerging, which means that a stock was not
exploited considerably before, or we otherwise nominate it as mature.  
If a region is determined in this way, an appropriate
control rule can be selected.

Whenever the system stays
within the viability domain only the rules (0, partly), (1) and
(3) come into play which warrants that the fishery
will remain in the viability domain. 
In contrast, if a fishery leaves the viability domain, the
rules (0), (2), (3), and (4) are able to steer the industry 
back to sustainable limits. 

In case of an emerging fishery, we start with a high $x$ (moving
from the right to the left in Fig. \ref{fig:phase12a}). 
In region (3) it could be reasonable to
increase $r$ in order to harvest as much as possible. 
However, there would be a significant risk to move close to region (0), possibly
without recognizing qualitatively that the stock
decreases below $\underline{x}$. 
Such a situation can be avoided if the
scientific institution strictly follows the conservative control scheme in an emerging
fishery until region (1) is reached for the first time. Thenceforward we
call a fishery mature and catch recommendations can be increased to
improve the profits. Once
the fish stock begins to decrease again, we know that we have approached 
region (3) and $r$ must reduced to force the fishery back to (1) in order to 
stay within the viability domain. Contrary to the emerging fishery, it
is not needed to decrease $r$ until $h=\underline{h}$, because the
fish stock begins to increase before (region 1).

Now suppose that we are outside the viability domain.
This may be a consequence of, for example, an ill-adapted management scheme
which was applied before the co-management came into play.
If the fishery stays in region (2) the scientific institution shall 
increase $r$ to obtain at least $\underline{h}$ as yield. If
decreasing stocks and low catches are observed in a mature fishery  
(region 4,) by considering the precautionary principle, we must assume that
$x < \underline{x}$. In such a case only a moratorium ($h \equiv 0$) can steer the
fishery back to sustainability. 
If $\bar{r}(x)$ has zeros and the mature fishery resides below the
smallest zero, an increase of $r$ results in entering region (4) (same
situation 
as above).  As an additional precondition for a safe development, 
it is required that recommendations are always binding: once $r>\hat{r}(x)$,
i.e. the fisheries harvest less than it was recommended $r$ shall be reduced.
This prevents that the stock size decreases below $a$ (cf. Fig.
\ref{fig:phase12}). If $b<a$ this guarantees that $x$ approaches the
viability domain. 

The qualitative control as presented here allows higher profits
than for conservative control, because recommendations are
increased above $\underline{r}(x)$ for max$(a,b)<x<c$, which results in
higher harvest rates. It is also more secure,
because it forces the system back to the viability
domain -- even if quantitative observations on the system are not
available.
\section{Discussion}
The previous analysis shows that participatory co-management schemes are not a priori
viable, since the outcome strongly depends on the relation between biological,
economic and political factors, and, in particular, on 
the catch recommendations of the scientific institution. In the future,
is is clear that ill-managed fisheries will radically reduce the self-determining
options of the coming generations in the fishing industry. The applied
viability concept shows how the dangerous effects related to measurement deficits
can be surmounted. In this way, the corresponding uncertainty can be 
confined to a minimum by a parallel achievement of management options.

In this contribution we have shown this for the three cases of ichtyocentric,
conservative, and qualitative control. 
An extreme case in this context is a recommendation strategy which is purely based on the
observation of fish stocks. This exposes the fishery to a high risk of economic
and ecological decline.  Such a situation can be substantially improved by
designing a more flexible strategy which only needs qualitative information
about the state of the fishery and does not deterministically fix the scientific
institution.   
In addition, it is shown that even for a fishery outside of 
a viable zone, there exists a good chance that it can be steered into the safe
region if a suitable control scheme is applied.   
A disadvantage of the presented model is that it does not take the costs of
change into account, i.e.  that rapid changes in harvest recommendations may
induce high adaption or political costs. On the other hand, even under
uncertainty the
qualitative control strategy is at least as good as economically conservative control, and less
risky than data-intensive ichthyocentric management.

However, the knowledge regarding relevant processes 
in specific domains of marine fisheries will remain less over the 
next few decades. Referring to
the variety of deficits in environmental management discussed here, 
the situation in marine capture fisheries more explicitly shows
that we need additional techniques to enhance our knowledge.
This holds, of course, if we deal with higher-dimensional
systems. Here qualitative simulation offers a possibility of 
revealing further insights into the phase space structure under 
uncertainty \shortcite{Kuipers.1994}. It takes
assumptions about the monotonicity of the right-hand sides of an ODE
and generates all temporal sequences of trends and thresholds which
are consistent with the assumptions. The result of such a simulation
is a directed graph, where each node represents a region of the phase space
(e.g. the regions (1) to (4) in figure \ref{fig:phase12a}), and each
edge indicates a possible trajectory crossing the boundary
between the associated regions (for recent examples in fisheries, cf. 
\shortciteNP{Kropp.2001,Kropp.2002a,Kropp.2005d}).
By investigating the subgraph induced by the nodes which represent
regions adjacent to the boundary of a viability domain, we can find out where
the system is on the safe side, where it necessarily degrades, and
where the result depends on the control strategy.
\section{Conclusions}
The situation in industrial fisheries is still facing
serious challenges in both an ecological and economic sense.
Management objectives are rarely achieved in practice
and the debate about adequate management strategies is still
ongoing. Since a unique solution is not expected, there still exists the need for some kind of integrated assessment.

In this paper we presented a new approach in the development
and assessment of co-management regimes. It is 
based upon the experience that for sustainable fishery management steering strategies should take diverse uncertainties into account.
This is addressed by
extracting some robust system properties, even from weak
information, and by giving an overview of the capabilities of viability theory 
in sustainability research.
It is shown in this paper that
the viability of a fishery strongly depends on
the catch recommendations of a scientific institution
participating in the co-management framework. 
The applied methodology makes it feasible to develop a
qualitative control strategy which requires only little information 
about the state of the fishery and is less risky than data-rich management schemes.
Of course, this variety of outcomes is
not a general objection against co-management. For the future
we will use enhanced models, 
including capital and investment dynamics of the fishing
firms. We think that such a strategy not only introduces
additional inertia and modified rules for the
negotiation process, but also more flexible steering instruments (e.g. effort control).
Summarizing, we think that such a analytical concept paves the road toward
detailed insights into what happens in marine capture fisheries.
We expect similar valuable clues for more complex models, leading
towards an integrated assessment of fisheries, 
including ecological, economic, and social issues.
\newline\newline\noindent{\bf Acknowledgements:}
The authors wish to thank the anonymous reviewers for their comments which helped
us to improve the paper.
\appendix
\section{Propositions and Proofs}
\begin{propo} 
\label{prop:vd}
The interval $I=[\underline{x},\infty]$ is a viability domain, if and only if
\begin{enumerate}
  \item[(i)] $\hat{r}(\underline{x}) \ge \underline{r}(\underline{x})$, and
  \item[(ii)] $\underline{h} \le R(\underline{x})$, and
 \item[(iii)] $\underline{r}(\underline{x}) \ge 0$ or $\bar{r}(\underline{x}) \ge 0$.
\end{enumerate}
\end{propo}
\begin{proof}
We start with some technical remarks: 
If $r \ge \underline{r}$ then
$h_b(x,r) \ge h_b(x, \underline{r}(x))$, since $D_r h_b < 0$.
The latter is equivalent to $\underline{h} > 0$ which can be shown
by basic
calculations, i.e.
\begin{equation}
  \label{eq:h1pos}
 h_b(x,r) > 0.
\end{equation}
In particular, it follows from condition (\ref{prop:vd}i) 
and Eq. (\ref{eq:hath1h2}) that
\begin{equation}
  \label{eq:hatrpos}
 \hat{r}(x) = h_b(x, \hat{r}(x)) > 0.
\end{equation}
\begin{enumerate}
\item[1.] It is shown that $I$ is a viability domain, if the 
conditions (\ref{prop:vd}i-iii) are valid.
As a consequence of the monotonicity of $\hat{r}(x)$ and
$\underline{r}(x)$ (Eq. \ref{eq:Dxrhat} and \ref{eq:Dxru}) it follows from condition
(i) and from $x\ge\underline{x}$ that also 
 $\hat{r}(\underline{x}) \ge \underline{r}(\underline{x})$.

We have to show that $\forall x > \underline{x} ~\exists\, r \ge 0:
h(x,r)\ge\underline{h}$ and that $\exists\, r\ge 0: h(\underline{x},r)
\ge \underline{h}$ and $R(\underline{x}) \ge h(\underline{x},r)$. We
claim that $r=\max(0,\underline{r}(x))$ respects these properties.

If $\underline{r}(\underline{x}) \ge 0$ (condition \ref{prop:vd}iii),
$r=\underline{r}(x)$. Then
$h(x,r) = h_b(x,r)$, due to Eq. (\ref{eq:harvest}),
$r \le \hat{r}(x)$, and Eq. (\ref{eq:h1pos}). By
elementary calculations we obtain $h_b(x,\underline{r}(x))=\underline{h}$. In
particular, for $x=\underline{x}$, we have
$h_b(x,\underline{r}(\underline{x}))=\underline{h}\le R(\underline{x})$, because of condition (\ref{prop:vd}ii).

If $\underline{r}(\underline{x}) < 0$, by condition (iii), 
$\bar{r}(\underline{x}) \ge 0$, and $r=0$. Then $h(x,r)=h_b(x,0)$ due
to Eq. (\ref{eq:harvest}), $0 < \hat{r}(x)$, and
Eq. (\ref{eq:h1pos}). Thus, for all $x \ge \underline{x}$: $h_b(x,0) \ge
h_b(x,\underline{r}(x))=\underline{h}$. Moreover,
$h_b(\underline{x},0) \le R(\underline{x})$, because this is
equivalent to $\bar{r}(\underline{x}) \ge 0$.

\item[2.] Conversely, conditions (\ref{prop:vd}i-iii) are a consequence of
$I$ being a viability domain, i.e. $\forall x \ge
\underline{x} ~\exists\, r\ge 0: h(x,r)\ge \underline{h} > 0$ and if
$x=\underline{x}$ additionally $h(\underline{x},r)\le R(\underline{x})$.

It is $h_b(\underline{x},\underline{r}(\underline{x})) = \underline{h}
\le h(x,r) \le \hat{r}(\underline{x}))$. The first inequality is valid by
assumption, the second due to Eq. (\ref{eq:hlehatr}). Considering
Eq. (\ref{eq:hath1h2}) the last term is equivalent to
$h_b(\underline{x},\hat{r}(\underline{x}))$. Due to a strict monotonicity of
$h_b$ in $r$ condition (\ref{prop:vd}i) follows. 

It follows directly from the assumptions that $\underline{h} \le
h(\underline{x},r) \le R(\underline{x})$, i.e. condition (ii).

Now let us suppose that $\underline{r}(\underline{x}) < 0$.
Equation (\ref{eq:hath1h2}) and condition (\ref{prop:vd}i) obviously show that
$\hat{r}(\underline{x}) = h_b(\underline{x},\hat{r}(\underline{x}))
\ge h_b(\underline{x},r)) = \underline{h} > 0$. 
If $r \le \hat{r}(\underline{x})$ then $h_b(\underline{x},0)) \le
h_b(\underline{x},r)) = h(\underline{x},r)) \le R(\underline{x})$.
If $r > \hat{r}(\underline{x})$ then $h_b(\underline{x},0)) \le
h_b(\underline{x},\hat{r}(\underline{x})) =
h_n(\underline{x},\hat{r}(\underline{x})) =
h_n(\underline{x},r)\le R(\underline{x})$. Thus, in any case 
$h_b(\underline{x},0)) \le  R(\underline{x})$, which is equivalent to
$\bar{r}(\underline{x}) \ge 0$, condition (\ref{prop:vd}iii) holds. \hfill $\blacksquare$
\end{enumerate}
\label{proof:vd}
\end{proof}
\begin{propo}
If $x \in I=[\underline{x}, \infty]$, $I$ a viability domain and consequently 
$\underline{r}(x) \le \hat{r}(x)$. $r$ is economically viable, if and
only if $R(x) \ge \underline{r}(x)$.
\label{propo:icht1}
\end{propo}
\begin{proof}
At first we
observe that $h_b(x,R(x)) \ge h_b(x,
\underline{r}(x))=\underline{h}>0$ and
$\hat{r}(x)=h_b(x,\hat{r}(x)) \ge h_b(x,\underline{r}(x))
=\underline{h}>0$.
Thus, if $r=R(x) \ge \hat{r}(x)$ it follows (with
Eq. \ref{eq:harvest}, \ref{eq:hath1h2} and \ref{eq:Drh1}) that
$h(x,R(x)) = h_n(x,R(x)) = \hat{r}(x) = h_b(x,\hat{r}(x)) \ge
 h_b(x,\underline{r}(x)) = \underline{h}$. Also, if $\underline{r}(x) \le
 R(x) \le \hat{r}(x)$ then $h(x,R(x))=h_b(x,R(x))$ $\ge
 h_b(x,\underline{r}(x)) = \underline{h}$. Contrary, if $R(x) \le
 \underline{r}(x)$, then $h(x,R(x))$ equals $0 < \underline{h}$ or
 $h_b(x,R(x)) \le h_b(x,\underline{r}(x)) = \underline{h}$. \hfill$\blacksquare$
\label{proof:icht1}
\end{proof} 
\begin{propo}
$h(x,R(x)) \le R(x)$, if and only if $r=R(x) \ge \hat{r}(x)$.
\label{propo:icht2}
\end{propo}
\begin{proof}

If $R(x) \ge \hat{r}(x)$, then $h(x,R(x))=h_n(x,R(x))=$
$\hat{r}(x) \le R(x)$. 
Otherwise, if $R(x) < \hat{r}(x)$
we use the fact that $\hat{r}(x) > r \Leftrightarrow h_b(x,r)>r$
(true by elementary calculations), to see that $h_b(x,R(x)) > R(x) \ge
0$, and we have $h(x,R(x)) =h_b(x,R(x)) > R(x)$ (by Eq. \ref{eq:harvest}). \hfill$\blacksquare$
\label{proof:icht2}
\end{proof} 

\end{document}